\documentclass[12pt]{article}

\usepackage{todonotes}

\usepackage{amsmath,amssymb,amsbsy,amsfonts,amsthm,latexsym,
                        amsopn,amstext,amsxtra,euscript,amscd,color}
\usepackage{color,xcolor}
\usepackage{latexsym}
\usepackage{cite}
\usepackage{amsfonts}
\usepackage{amsfonts,mathrsfs}
\usepackage{graphicx}
\usepackage{psfrag}
\usepackage{subfigure}
\usepackage{url}
\usepackage{stfloats}
\usepackage{amsmath}
\usepackage{algorithm}
\usepackage{algorithmic}
\usepackage{hyperref}

\hypersetup{colorlinks=true}

\newtheorem{theorem}{Theorem}
\newtheorem{lemma}{Lemma}

\newcommand{\quash}[1]{}

\begin{document}

\title{Probabilistic results on the $2$-adic complexity}

\author{Zhixiong~Chen~$^1$ and Arne~Winterhof~$^2$\\
~\\
$^1$ Fujian Key Laboratory of Financial Information Processing,\\
 Putian University, Putian, Fujian 351100, P. R. China\\
$^2$ Johann Radon Institute for Computational and Applied Mathematics,\\
Austrian Academy of Sciences, Altenbergestr.\ 69,
A-4040 Linz, Austria
}

\maketitle

\begin{abstract}

This work is devoted to solving some closely related open problems on the average and asymptotic behavior of the $2$-adic complexity of binary sequences.
First, for fixed $N$, we prove that  the expected value~$E^{\mathrm{2-adic}}_N$ of the $2$-adic complexity over all binary sequences of length $N$ is close to $\frac{N}{2}$ and the deviation from $\frac{N}{2}$ is at most of order of magnitude $\log(N)$. More precisely, we show that
  $$\frac{N}{2}-1 \le E^{\mathrm{2-adic}}_N= \frac{N}{2}+O(\log(N)).$$
  We %will
  also prove bounds on the expected value of the $N$th rational complexity.

Our second contribution is to prove
for a random binary sequence $\mathcal{S}$ that the $N$th $2$-adic complexity satisfies with probability~1
$$
\lambda_{\mathcal{S}}(N)=\frac{N}{2}+O(\log(N)) \quad \mbox{for all $N$}.
$$



\end{abstract}

\section{Introduction}

Binary sequences
are used for many cryptographic applications
such as pseudorandom number generators for stream cipher cryptosystems. An output sequence, denoted by $\mathcal{S}=(s_i)_{i\geq 0}$, of a linear feedback shift register (LFSR) of length $L$ is
 determined by a recurrence relation with $a_j\in\{0,1\}$ for $0\leq j\leq L-1$,
$$
s_{i+L}=a_{L-1}s_{i+L-1}+a_{L-2}s_{i+L-2}+\ldots+a_0s_{i} \bmod{2},\quad i\geq 0.
$$
Sequences derived from a short LFSR are considered weak in view of cryptographic applications. In order to determine this shortness, the notion of \textit{Nth linear complexity} of $\mathcal{S}$, denoted by $~L_{\mathcal{S}^N}(N)$ was introduced, see for example Rueppel \cite{R85,R86},
as the length $L$ of a shortest LFSR that generates $\mathcal{S}^{N}=(s_0,s_1,\ldots,s_{N-1})$, the first
$N$ elements of $\mathcal{S}$. Sometimes we will use the notation $L_{\mathcal{S}}(N)=L_{\mathcal{S}^N}(N)$.

Rueppel \cite{R86}, see also \cite[Sect.\ 18.1.1]{GK2012}, proved that, for fixed $N$, the expected value of the $N$th linear complexity
$$E_N^{\mathrm{lin}}=\frac{1}{2^N}\sum_{L=0}^N  A_N^{\mathrm{lin}}(L)L$$
of binary sequences of length $N$
satisfies
\begin{equation}\label{lin1}
E_N^{\mathrm{lin}}=\frac{N}{2}+O(1),
\end{equation}
where $A_N^{\mathrm{lin}}(L)$ is the number of sequences $\mathcal{S}^N\in \{0,1\}^N$ with $L_{\mathcal{S}^N}(N)=L$. The result is based on the exact formula of the counting function
\begin{equation}\label{count}A_N^{\mathrm{lin}}(L)=2^{\min\{2L-1,2N-2L\}},\quad L=0,1,\ldots,N.
\end{equation}
Niederreiter \cite{N1988,N1990} proved that, for varying $N$, with probability $1$
\begin{equation}\label{lin2}
L_{\mathcal{S}}(N)=\frac{N}{2}+O(\log(N)) \quad \mathrm{for}~~ N\geq 2
\end{equation}
for a random sequence $\mathcal{S}$.\\

Goresky and Klapper proposed a new type of registers called feedback with carry shift registers (FCSRs) in the 1990s, see \cite{GK2012,KG1997}. The algorithm (of an FCSR of length~$L$) is carried out by  a recurrence relation in the ring $\mathbb{Z}$ of integers with $a_j\in\{0,1\}$ for $0\leq j\leq L-1$,
$$
s_{i+L}+2z_{i+L}=a_{L-1}s_{i+L-1}+a_{L-2}s_{i+L-2}+\ldots+a_0s_{i}+z_{i+L-1}, \quad i\geq 0,
$$
with
initial integer $z_{L-1}$ and carries $z_{L},z_{L+1},\ldots$
The \textit{$N$th $2$-adic complexity} $\lambda_{\mathcal{S}^N}(N)$ of $\mathcal{S}$
characterizes the length of a shortest FCSR that generates the  first $N$ elements of $\mathcal{S}$.
More precisely, $\lambda_{\mathcal{S}}(N)=\lambda_{\mathcal{S}^N}(N)$ is defined as the binary logarithm of the
{\em $N$th rational complexity}
$$\Lambda_{\mathcal{S}}(N)=
\Lambda_{\mathcal{S}^N}(N)=\min \left\{\max\{|f|, q\} : f,q\in\mathbb{Z}, q>0 \mbox{ odd}, q\sum_{i=0}^{N-1}s_i2^i \equiv f \bmod{2^N}   \right\};
$$
see \cite[p.\ 328]{GK2012}, the recent survey \cite{W2023} or \cite{TQ2010}, that is,
$$\Lambda_{\mathcal{S}^N}(N)=2^{\lambda_{\mathcal{S}^N}(N)}.$$
(Note that, more recently, the term rational complexity has been introduced for a different measure of pseudorandomness, see \cite{VCJ2022}. More precisely, \cite{VCJ2022} identifies binary sequences with reals in the unit interval whereas here we identify binary sequences with $2$-adic integers.)
Any pair $(q,f)$
with $\Lambda_{\mathcal{S}^N}(N)=\max\{|f|,q\}$ is called a \textit{minimal rational representation} for the first $N$ elements of $\mathcal{S}$.
It is trivial that
$$\lambda_{\mathcal{S}}(N)\leq \lambda_{\mathcal{S}}(N+1)\quad  \mbox{for $N\geq 1$}.$$

Although experimental results support the conjecture that the analog of $(\ref{lin1})$ for the expected value of the $N$th $2$-adic complexity
holds true, this problem, stated in \cite{GK2012,K2014,TQ2010}, is still open. Some partial results can be found in \cite{K2007}, \cite[Sect.\ 18.5]{GK2012} and \cite{TQ2010}. In particular, Tian and Qi \cite{TQ2010} proved bounds on the expected value
$$E_N^{\mathrm{rat}}=\frac{1}{2^N}\sum_{\mathcal{S}^N\in \{0,1\}^N}\Lambda_{\mathcal{S}^N}(N)$$
of the $N$th rational complexity,
\begin{equation}\label{TQbound}2^{0.7716N}>E_N^{\mathrm{rat}}\ge \left\{\begin{array}{cc}\frac{2^{(N+1)/2}}{9\log_2(N-1)-1}+\frac{N-1}{4},& N \mbox{ odd},\\
\frac{2^{(N-2)/2}}{9\log_2(N-2)-1}+\frac{N}{4},& N \mbox{ even},
\end{array} \right.\quad N\ge 3,\end{equation}
which led to an upper bound on the expected value
$$E_N^{2-\mathrm{adic}}
=\frac{1}{2^{N}}\sum\limits_{\mathcal{S}^N\in\{0,1\}^N}\lambda_{\mathcal{S}^N}(N)$$
of the $N$th $2$-adic complexity,
\begin{equation}\label{TQbound2}
E_N^{2-\mathrm{adic}}
<0.7716N.
\end{equation}
However, their experiments indicated that  $E_N^{2-\mathrm{adic}}$ is near $\frac{N-1}{2}$.
Unfortunately, an analog of (\ref{count}) for the $2$-adic complexity seems to be out of reach and we have to deal instead with estimates on the number of sequences with `small' and `large' $2$-adic complexity, respectively, which still leads to the following substantial improvement.

\begin{theorem}\label{expected-value}
The expected values $E_N^{\mathrm{rat}}$ and $E^{\mathrm{2-adic}}_N$ over all binary sequences of length~$N$ satisfy
$$E_N^{\mathrm{rat}}= 2^{N/2+O(N/\log(N))}$$
and
$$
E^{\mathrm{2-adic}}_N=\frac{N}{2}+O(\log(N)).$$ 
\end{theorem}

The problem of finding an analog of $(\ref{lin2})$, stated in \cite[p.392]{GK2012} and \cite[Sect.8]{W2023}, is also settled in this work. That is, we prove
for a random binary sequence $\mathcal{S}$ that %the $N$th $2$-adic complexity
$\lambda_{\mathcal{S}}(N)$ is close to $\frac{N}{2}$
and for all $N$ the deviation from $\frac{N}{2}$
is at most of order of magnitude $\log(N)$ with probability~1. This also implies
$$\Lambda_{\mathcal{S}}(N)=2^{N/2+O(\log(N))}\quad \mbox{for all } N \mbox{ with probability } 1,$$
see Theorem \ref{Thm-asymp} in Section~\ref{Asy-2adic}.
Vielhaber \cite[Corollary~21]{V2007} considered the \textit{$2$–adic jump complexity} of $\mathcal{S}^{N}$, the
number of changes of the minimal rational representation $(q_k,f_k)$ for $\mathcal{S}^{k}$ for all $0\le k<N$.

In the sequel, as in \cite{N1990}, we will make use  of
the equiprobable space of all finite sequences of the
same length over $\{0,1\}$.  In detail, let $\mu$ be the uniform probability measure on $\{0,1\}$ which assigns the measure $\frac{1}{2}$ to each element of $\{0,1\}$. This measure
induces the complete product measure $\mu^{\infty}$
on the set $\{0,1\}^{\infty}$ of infinite sequences over~$\{0,1\}$.

Moreover, we will use the notation
$$f(N)=O(g(N))\quad \mbox{if}\quad |f(N)|\le cg(N)\quad  \mbox{for some constant $c>0$}$$
and
$$f(N)=o(g(N))\quad \mbox{if} \quad \lim\limits_{N\rightarrow \infty} \frac{f(N)}{g(N)}=0.$$
Sometimes we also use $f(N)\ll g(N)$ and $g(N)\gg f(N)$ instead of $f(N)=O(g(N))$.

In Section~\ref{small_and_large} we estimate the number of sequences with `small' and `large' $2$-adic complexity, respectively. Based on these bounds we derive in Section~\ref{expect} a lower and an upper bound on the expected value, from which we can easily derive Theorem~\ref{expected-value}.
In Section~\ref{Asy-2adic} we study the asymptotic behavior of the $N$th $2$-adic complexity for random binary sequences.

\section{Numbers of sequences with certain restrictions}
\label{small_and_large}

For a positive integer $w$, let
$M_N(w)$ be the number of $\mathcal{S}^N\in \{0,1\}^N$ with $N$th $2$-adic complexity
$\log_2(w)$ (or $N$th rational complexity $w$),
\begin{eqnarray*}
M_N(w)&=&\left|\left\{\mathcal{S}^N\in \{0,1\}^N: \lambda_{\mathcal{S}^N}(N)=\log_2(w) \right\}\right|\\
&=&\left|\left\{\mathcal{S}^N\in \{0,1\}^N: \Lambda_{\mathcal{S}^N}(N)=w \right\}\right|.
\end{eqnarray*}

First, we estimate the number of $\mathcal{S}^N$ with `small' $N$th $2$-adic complexity, that is, $\lambda_{\mathcal{S}^N}(N)<\frac{N}{2}-\delta\log_2(N)$.

\begin{lemma}\label{num-lowerbound}
We have
$$\sum_{w=1}^W M_N(w)=\frac{8}{\pi^2}W^2+O(W\log(W)),\quad W\le 2^{(N-1)/2},$$
where the implied constant is absolute.

In particular,
for  any $\delta>0$ and
$$
\Delta_{N}=\left\{\mathcal{S}^N\in \{0,1\}^{N}: \lambda_{\mathcal{S}^N}(N)<\frac{N}{2}-\delta\log_2(N) \right\},
$$
we have
$$
|\Delta_{N}|=
 \frac{8}{\pi^2}\frac{2^N}{N^{2\delta}}+O\left(\frac{2^{\frac{N}{2}}}{N^{\delta-1}}\right),
$$
where the implied constant is absolute.
\end{lemma}

 Proof.
 For $w\le 2^{(N-1)/2}$ we can express $M_N(w)$ in terms of Euler's totient function
$$M_N(w)=\left\{\begin{array}{cc}
2\varphi(w), & w \mbox{ even},\\
3\varphi(w), & w \mbox{ odd},
\end{array}\right. \quad w\le 2^{(N-1)/2},$$
where
$$\varphi(w)=|\{x=1,2,\ldots,w: \gcd(x,w)=1\}|,\quad w=1,2,\ldots,$$
see \cite[Lemmas 2 and 3]{TQ2010}.
Hence, we have
$$\sum_{w=1}^WM_N(w)=3\sum_{w=1}^W \varphi(w)-\sum_{w=1}^{\lfloor W/2\rfloor}\varphi(2w),\quad W\le 2^{(N-1)/2}-1.$$
Since
$$\varphi(2w)=\left\{\begin{array}{cc} \varphi(w), & w \mbox{ odd},\\
2\varphi(w), & w \mbox{ even},
\end{array}\right.$$
we get
\begin{eqnarray*}\sum_{w=1}^{\lfloor W/2\rfloor}\varphi(2w)&=&
\sum_{w=1}^{\lfloor W/2\rfloor}\varphi(w)+\sum_{w=1}^{\lfloor W/4\rfloor}\varphi(2w)\\
&=&\sum_{w=1}^{\lfloor W/2\rfloor}\varphi(w)+\sum_{w=1}^{\lfloor W/4\rfloor}\varphi(w)+\sum_{w=1}^{\lfloor W/8\rfloor}\varphi(2w)\\
&=& \cdots\\
&=&\sum_{j=1}^{\lfloor \log_2(W)\rfloor}\sum_{w=1}^{\lfloor W/2^j\rfloor}\varphi(w).
\end{eqnarray*}

Using Walfisz' formula, see for example \cite[Theorem~330]{HW1979}, the original work \cite{W1963}, or \cite[Theorem~1]{L2016} for an improved error term,
$$\sum\limits_{w=1}^{m}\varphi(w) =\frac{3}{\pi^2}m^2+O\left(m \log(m) \right),$$
 we get
\begin{eqnarray*}
\sum_{w=1}^WM_N(w)
=  RW^2+O(W\log(W)),\quad W\le 2^{(N-1)/2},
\end{eqnarray*}
where
$$R=\frac{3}{\pi^2}\left(3-\sum_{j=1}^{\lfloor \log_2(W)\rfloor}4^{-j}\right)=\frac{8}{\pi^2}+O(W^{-2})$$
and the result follows. \qed \\

Now, we estimate the number of $\mathcal{S}^N$ with `large' $N$th $2$-adic complexity, that is, $\lambda_{\mathcal{S}^N}(N)>\frac{N}{2}+\delta\log_2(N)$. We will need the following result on the growth of $\Lambda_\mathcal{S}(N)$,   see \cite[Lemma 18.5.1]{GK2012}.

\begin{lemma}\label{star}
We have
$$
\Lambda_{\mathcal{S}}(N)\leq \Lambda_{\mathcal{S}}(N-1)+ \frac{2^{N-1}}{\Lambda_{\mathcal{S}}(N-1)},\quad N\ge 2.
$$
\end{lemma}

A pivotal fact is that we will use the solutions of the inequality in Lemma~\ref{star} in the proof of Lemma~\ref{num-upperbound}, below, that is, the
 solutions $x>0$ of an inequality of the form
 $$x+\frac{a}{x}\geq b\quad \mbox{or}\quad x^2-bx+a\geq 0,$$
 where $0<4a\le b^2$, namely,
$$0<x\leq \frac{b-\sqrt{b^2-4a}}{2} ~~~\mathrm{or}~~~ x\geq \frac{b+\sqrt{b^2-4a}}{2}.$$

We will need another preliminary result on sums of reciprocals.
\begin{lemma}\label{func-N}
For $\delta>\frac{1}{2}$
and $0< \varepsilon<1$ we have
$$
\sum\limits_{c=0}^{\lfloor\varepsilon N\rfloor}\frac{1}{(N-c)^{2\delta}}\le  \frac{c(\varepsilon,\delta)}{N^{2\delta-1}},
$$
where
$$c(\varepsilon,\delta)=\frac{(1-\varepsilon)^{1-2\delta}-1}{2\delta-1}=2^{O(\delta)}.$$
\end{lemma}
Proof. We get the result from
\begin{eqnarray*}
\sum\limits_{c=0}^{\lfloor\varepsilon N\rfloor}\frac{1}{(N-c)^{2\delta}}&=& \sum\limits_{x=N-\lfloor\varepsilon N\rfloor}^{N}\frac{1}{x^{2\delta}}\\
&\leq & \int^{N}_{N-\varepsilon N} x^{-2\delta} \mathrm{d}x\\
&=& -\frac{1}{2\delta-1} \cdot \left(N^{1-2\delta}-(N-\varepsilon N)^{1-2\delta} \right)\\
&=& c(\varepsilon,\delta) N^{1-2\delta}.
\end{eqnarray*}
Note that the constant $c(\varepsilon,\delta)$ is positive under the conditions on $\delta$ and $\varepsilon$.
\qed\\

Now we are ready to prove an upper bound on the number of sequences with
`large' $2$-adic complexity.
\begin{lemma}\label{num-upperbound}
Consider
$$
\Gamma_{N}=\left\{\mathcal{S}^N\in \{0,1\}^{N}: \lambda_{\mathcal{S}^N}(N)>\frac{N}{2}+\delta\log_2(N) \right\}.
$$
For $\delta >\frac{1}{2}$
we have
$$
|\Gamma_{N}| \le \frac{2^{N+O(\delta)}}{N^{2\delta-1}}.
$$
where the implied constant is absolute.
\end{lemma}
 Proof. For any integer $c\ge 0$ we consider the number of sequences $\mathcal{S}^{N-c}$, the first $N-c$ elements of $\mathcal{S}^N$, with `large' $\lambda_{\mathcal{S}^{N-c}}(N-c)$,
$$
K_N(c)=\left|\left\{\mathcal{S}^{N-c}\in \{0,1\}^{N-c}: \lambda_{\mathcal{S}^{N-c}}(N-c)>\frac{N}{2}+\delta\log_2(N)-0.4c \right\}\right|,
$$
that is, $|\Gamma_N|=K_N(0)$.
Since
$\lambda_{\mathcal{S}^{N-c}}(N-c)\le N-c-1$
by \cite[Lemma 10]{TQ2010}
we have
$$
K_N(c)=0\quad \mbox{for}\quad c\ge \frac{5}{6}N-1.
$$
Hence, below we restrict ourselves to the case $c<\frac{5}{6}N-1<N-3$, that is, $N>12$.

Now from
$$\lambda_{\mathcal{S}^{N-c}}(N-c)>\frac{N}{2}+\delta\log_2(N)-0.4c,$$
we get
by Lemma \ref{star},
\begin{eqnarray*}
2^{\frac{N-c}{2}+\delta\log_2(N-c)}&<&2^{\frac{N}{2}+\delta\log_2(N)-0.4c}< 2^{\lambda_{\mathcal{S}^{N-c}}(N-c)}\\
&\leq&  \Lambda_{\mathcal{S}^{N-c-1}}(N-c-1)+ \frac{2^{N-c-1}}{\Lambda_{\mathcal{S}^{N-c-1}}(N-c-1)},
\end{eqnarray*}
which implies either
\begin{eqnarray*}
\Lambda_{\mathcal{S}^{N-c-1}}(N-c-1)&<&2^{\frac{N-c}{2}+\delta\log_2(N-c)-1}-\sqrt{2^{N-c+2\delta\log_2(N-c)-2}-2^{N-c-1}}\\
&=&2^{\frac{N-c}{2}+\delta\log_2(N-c)-1}\left(1-\sqrt{1-2^{1-2\delta\log_2(N-c)}}\right)\\
&=&\frac{2^{\frac{N-c}{2}+\delta\log_2(N-c)-1}\cdot 2^{1-2\delta\log_2(N-c)}}{1+\sqrt{1-2^{1-2\delta\log_2(N-c)}}}\\
&<&\frac{2^{\frac{N-c}{2}-\delta\log_2(N-c)}}{\sqrt{2}},\quad (\mathrm{for}~~ c\le N-3),\\
&<& 2^{\frac{N-c-1}{2}-\delta\log_2(N-c-1)}
\end{eqnarray*}
or
\begin{eqnarray*}
\Lambda_{\mathcal{S}^{N-c-1}}(N-c-1)&>&2^{\frac{N}{2}+\delta\log_2(N)-0.4c-1}+\sqrt{2^{N+2\delta\log_2(N)-0.8c-2}-2^{N-c-1}}\\
& = & 2^{\frac{N}{2}+\delta\log_2(N)-0.4c-1}\left(1+\sqrt{1-2^{1-2\delta\log_2(N)-0.2c}}\right)\\
 & \geq & 2^{\frac{N}{2}+\delta\log_2(N)-1-0.4c}\cdot 2^{1-0.4},\quad (\mathrm{for}
 ~~ N\ge 3),\\
& = & 2^{\frac{N}{2}+\delta\log_2(N)-0.4(c+1)}.
\end{eqnarray*}

Now for a fixed positive integer $C\le \left\lfloor\frac{5}{6}N\right\rfloor$ assume that
$$\lambda_{\mathcal{S}^{N-c}}(N-c)>\frac{N}{2}+\delta\log_2(N)-0.4c\quad \mbox{for }c=0,1,\ldots,C-1$$
and
\begin{equation}\label{small}
\lambda_{\mathcal{S}^{N-c}}(N-C)<2^{\frac{N-C}{2}-\delta \log_2(N-C)}.
\end{equation}
There are $O\left(\frac{2^{N-C}} {(N-C)^{2\delta}}\right)$
sequences $\mathcal{S}^{N-C}$ with `small' $\lambda_{\mathcal{S}^{N-C}}(N-C)$ satisfying $(\ref{small})$ by Lemma~\ref{num-lowerbound}. Each such $\mathcal{S}^{N-C}$ can be padded to a sequence $\mathcal{S}^N$ of length $N$ in $2^C$ ways and  there are $O\left(\frac{2^N}{(N-C)^{2\delta}}\right)$ sequences with `large' $\lambda_{\mathcal{S}^N}(N)$.
Hence,
\begin{eqnarray*}
|\Gamma_N| \ll  2^N\sum\limits_{C=1}^{\left\lfloor \frac{5}{6}N\right\rfloor}
\frac{1}{(N-C)^{2\delta}}.
\end{eqnarray*}
Finally, Lemma~\ref{func-N} (with $\varepsilon=\frac{5}{6})$ implies the desired result.  \qed

\section{Expected values of $2$-adic and rational complexity}
\label{expect}

In this section, we prove lower and upper bounds on the expected values of the $N$th $2$-adic and rational complexity. Their combination implies Theorem~\ref{expected-value}.
We first prove lower bounds  based on the following lemma.

\begin{lemma}\label{Le-111}
For $N\geq 2$ and any $\mathcal{S}^{N-1}=(s_0,s_1,\ldots,s_{N-2})\in \{0,1\}^{N-1}$,  we have
$$\lambda_{(\mathcal{S}^{N-1},0)}(N)+\lambda_{(\mathcal{S}^{N-1},1)}(N)\ge
N-2$$
 and
$$\Lambda_{(\mathcal{S}^{N-1},0)}(N)+\Lambda_{(\mathcal{S}^{N-1},1)}(N)\ge
2^{\frac{N}{2}},$$
where $(\mathcal{S}^{N-1},a)=(s_0,s_1,\ldots,s_{N-2},a)$ for $a\in \{0,1\}$.
\end{lemma}
Proof. It suffices to show $$\Lambda_{(\mathcal{S}^{N-1},0)}(N)\cdot\Lambda_{(\mathcal{S}^{N-1},1)}(N)\geq 2^{N-2},$$
from which the first inequality follows by
taking logarithms
and the second inequality follows by the inequality of arithmetic and geometric means
$$
\Lambda_{(\mathcal{S}^{N-1},0)}(N)+\Lambda_{(\mathcal{S}^{N-1},1)}(N)
\geq
2\sqrt{\Lambda_{(\mathcal{S}^{N-1},0)}(N)\cdot\Lambda_{(\mathcal{S}^{N-1},1)}(N)}.
$$

Now take minimal rational representations $(q,f)$ of $(\mathcal{S}^{N-1},0)$ and $(q',f')$ of $(\mathcal{S}^{N-1},1)$. Then we have
 $$q\sum_{i=0}^{N-2}s_i2^i\equiv f \bmod 2^{N} ~~ \mbox{ and }~~  q'\sum_{i=0}^{N-2}s_i2^i+2^{N-1}\equiv f' \bmod 2^{N}.$$
It follows that $fq'-f'q \equiv 2^{N-1}  \bmod{2^N}$
and hence
$$
 2^{N-1}\leq |fq'-f'q |\leq |fq'|+|f'q| \leq 2\Lambda_{(\mathcal{S}^{N-1},0)}(N)\cdot\Lambda_{(\mathcal{S}^{N-1},1)}(N),
$$
which completes the proof.  \qed\\

From Lemma \ref{Le-111} one can directly get
lower bounds on the expected values $E^{\mathrm{2-adic}}_N$ and $E^{\mathrm{rat}}_N$.

\begin{theorem}\label{ex-lowerbound}
We have
$$
E^{\mathrm{2-adic}}_N\geq \frac{N}{2}-1
$$
and
$$
E^{\mathrm{rat}}_N\geq 2^{\frac{N}{2}-1}+\frac{N-5}{4}.
$$
\end{theorem}
Proof.
Lemma \ref{Le-111} leads to
\begin{eqnarray*}
2^N \cdot E^{\mathrm{2-adic}}_N = \sum\limits_{\mathcal{S}^{N-1}\in\{0,1\}^{N-1}}\left(\lambda_{(\mathcal{S}^{N-1},0)}(N)+\lambda_{(\mathcal{S}^{N-1},1)}(N)\right)\geq (N-2)2^{N-1}
\end{eqnarray*}
and
$$
2^N \cdot E^{\mathrm{rat}}_N = \sum\limits_{\mathcal{S}^{N-1}\in\{0,1\}^{N-1}}\left(\Lambda_{(\mathcal{S}^{N-1},0)}(N)+\Lambda_{(\mathcal{S}^{N-1},1)}(N)\right)\geq 2^{\frac{N}{2}}\cdot 2^{N-1},
$$
which imply the desired lower bounds on $E^{\mathrm{2-adic}}_N$ as well as
\begin{equation}\label{exprat}
E^{\mathrm{rat}}_N\ge 2^{\frac{N}{2}-1}.
\end{equation}

For $E^{\mathrm{rat}}_N$,
we slightly improve this bound to
$$E_N^{\mathrm{rat}}\ge 2^{\frac{N}{2}-1}+\frac{N-5}{4}.$$
Let $\mathcal{S}_k^N=(s_0,s_1,\ldots,s_{N-1})$ be any sequence satisfying $s_n=0$ for $n=0,1,\ldots,k-1$ and $s_k=1$.
By \cite[Lemma 4]{TQ2010} we have for any $\frac{N}{2}\le k\le N-2$,
$$\Lambda_{\mathcal{S}_k^{N}}(N)=2^k.$$
There are $2^{N-k-1}$ such sequences. Their contribution to the expected value $E_N^{\mathrm{rat}}$
is
$$\frac{1}{2^N}\sum_{k=\left\lceil \frac{N}{2}\right\rceil}^{N-2}2^k\cdot 2^{N-k-1}
=\frac{1}{2^N}\left(N-1-\left\lceil \frac{N}{2}\right\rceil\right)2^{N-1}\ge \frac{N-3}{4}.$$
However, their contribution to the lower bound (\ref{exprat}) was exactly $2^{\frac{N}{2}}$ for each pair sequences $\mathcal{S}_k^N$with different $s_{N-1}$ and there are $2^{N-k-2}$ many such pairs, that is totally
$$\frac{1}{2^N} \sum_{k=\left\lceil \frac{N}{2}\right\rceil}^{N-2}2^{\frac{N}{2}}\cdot 2^{N-k-2}
=\frac{1}{2^{N/2}}\left(2^{N-\lceil N/2\rceil-1}-1\right)< \frac{1}{2}.$$
This completes the proof of the second result. \qed\\

The lower bound on the expected value of the $N$th $2$-adic complexity in Theorem~\ref{ex-lowerbound} is very close to
the conjectured lower bound $\frac{N-1}{2}$ supported by the data for $N\le 20$ given in
\cite[Fig.1 and Table 1]{TQ2010}.
The lower bound on the $N$th rational complexity improves  the lower bound of \cite{TQ2010} in (\ref{TQbound}) substantially.

Now we turn to prove upper bounds which are much better than the upper bounds of \cite{TQ2010}, see  $(\ref{TQbound})$ and $(\ref{TQbound2})$.

\begin{theorem} 
For $\delta>\frac{1}{2}$ we have
$$
E^{\mathrm{2-adic}}_N \leq \frac{N}{2}+\delta\log_2(N)+O(N^{2-2\delta}),
$$
where the implied constant depends only on $\delta$,
and
$$E_N^{\mathrm{rat}}\le 2^{N/2+O(N/\log(N))},$$
where the implied constant is absolute.
\end{theorem}
Proof.
Take $\Gamma_N$ from Lemma~\ref{num-upperbound}. Since $\lambda_{\mathcal{S}^N}(N)\leq N-1$ by \cite[Lemma 10]{TQ2010} and $\delta>\frac{1}{2}$, we have
$$
\frac{1}{2^N}\sum\limits_{\mathcal{S}^N\in \Gamma_N} \lambda_{\mathcal{S}^N}(N)\leq   \frac{(N-1)\cdot|\Gamma_N|}{2^N} \ll
\frac{1}{N^{2\delta-2}}.%=o(1).
$$
Moreover, we have
$$
\frac{1}{2^N}\sum\limits_{\mathcal{S}^N\in \{0,1\}^{N}\setminus \Gamma_N}  \lambda_{\mathcal{S}^N}(N)
\leq \frac{2^N-|\Gamma_N|}{2^N}\cdot \left(\frac{N}{2}+\delta\log_2(N)\right)<\frac{N}{2}+\delta\log_2(N)
$$
and the first result follows.

Now we turn to prove the bound on $E_N^{\mathrm{rat}}$. Taking
$\delta_1=\delta,\delta_2,\ldots,\delta_m$ with
$$\frac{1}{2}<\delta=\delta_1<\delta_2<\ldots<\delta_m=O\left(\frac{N}{\log(N)}\right),$$ we subdivide $\Gamma_N$ into $m$ subsets
$$
\Gamma^{(i)}_N=\left\{\mathcal{S}^N\in \{0,1\}^{N}: \frac{N}{2}+\delta_i\log_2(N) <\lambda_{\mathcal{S}^N}(N)\leq \frac{N}{2}+\delta_{i+1}\log_2(N) \right\}
$$
for $1\le i\leq m-1$ and
$$
\Gamma^{(m)}_{N}=\left\{\mathcal{S}^N\in \{0,1\}^{N}: \lambda_{\mathcal{S}^N}(N)>\frac{N}{2}+\delta_m\log_2(N) \right\}.
$$
Then by Lemma \ref{num-upperbound} we see that
$$
|\Gamma^{(i)}_N|\le \frac{2^{N+O(\delta_i)}}{N^{2\delta_i-1}}=2^{N-2\delta_i\log_2(N)+O(N/\log(N))}, ~~~1\le i\leq m.
$$
Hence we get
\begin{eqnarray*}
E_N^{\mathrm{rat}}&=& \frac{1}{2^N}\left(\sum\limits_{i=1}^{m}\sum\limits_{\mathcal{S}^N\in \Gamma^{(i)}_N} \Lambda_{\mathcal{S}^N}(N)+\sum\limits_{\mathcal{S}^N\in \{0,1\}^{N}\setminus \Gamma_N}  \Lambda_{\mathcal{S}^N}(N)\right) \\
&<& \frac{1}{2^N}\left(\sum\limits_{i=1}^{m-1}|\Gamma^{(i)}_N|\cdot 2^{\frac{N}{2}+\delta_{i+1}\log_2(N)}+|\Gamma^{(m)}_N|\cdot 2^{N}+2^N\cdot 2^{\frac{N}{2}+\delta\log_2(N)}\right)\\
&\le& 2^{\frac{N}{2}+\delta\log_2(N)}+ \sum\limits_{i=1}^{m-1} 2^{\frac{N}{2}+(\delta_{i+1}-2\delta_i)\log_2(N)+O(N/\log(N))}\\
&&
+2^{N-2\delta_m\log_2(N)+O(N/\log(N))}.
\end{eqnarray*}

Now we balance the first $m$ summands by choosing
$$\delta_{i+1}=2\delta_i+\delta,\quad i=1,2,\ldots,m-1,$$
that is
$$\delta_i=(2^i-1)\delta,\quad i=1,2,\ldots,m.$$
Then we balance the first and the last summand, up to the $O$-term,
choosing
$$
\delta=\frac{N}{(2^{m+2}-2)\log_2(N)}.
$$
Hence with $m=\lceil\log_2(\log_2(N)) \rceil$ we derive
$$
E_N^{\mathrm{rat}}\leq (m+1)2^{\frac{N}{2}+\frac{N}{2^{m+2}-2}+O(N/\log(N))}= 2^{\frac{N}{2}+O(N/\log(N))},
$$
which completes the proof. \qed

\section{Asymptotic behavior of $2$-adic and rational complexity}\label{Asy-2adic}

As for the proof of $(\ref{lin2})$ in \cite{N1990} for the $N$th linear complexity, the Borel-Cantelli Lemma, see for example \cite[p.\ 228]{Lo1963}, also plays an important role in studying the behavior of the $N$th $2$-adic complexity. We state it here for the convenience of the reader.

\begin{lemma}[Borel-Cantelli Lemma]\label{BCLemma}
If the sum of the probabilities %$\mathrm{Pr}(-)$
of the events $\{\Theta_N : N=1,2,\ldots\}$ is finite
$$
\sum\limits_{N=1}^{\infty}\mathrm{Pr}(\Theta_N)<\infty,
$$
then the probability that infinitely many of them occur is $0$, that is,
$$
\mathrm{Pr}\left(\limsup\limits_{N\rightarrow \infty}\Theta_N\right)=0.
$$
\end{lemma}

Then the asymptotic behavior of the $2$-adic complexity and the rational complexity of a random sequence $\mathcal{S}$ is described in the following theorem.
\begin{theorem}\label{Thm-asymp}
For a random sequence $\mathcal{S}$ we have with probability $1$,
$$
\lambda_{\mathcal{S}}(N)=\frac{N}{2}+O(\log(N))
$$
and
$$\Lambda_{\mathcal{S}}(N)=2^{\frac{N}{2}+O(\log(N))}$$
for all $N$.
\end{theorem}
Proof. Fix any $\delta>1$, take $\Delta_N$ from Lemma \ref{num-lowerbound} and $\Gamma_N$ from Lemma  \ref{num-upperbound}. For
 $\Theta_N=\Delta_N \cup\Gamma_N$ we get
$$
\mu^{\infty}(\Theta_N)=\frac{|\Delta_N|+|\Gamma_N|}{2^{N}}\ll \frac{1}{N^{2\delta}}+\frac{1}{N^{2\delta-1}} \ll \frac{1}{N^{2\delta-1}}.
$$
Since $\delta>1$ we see that
$$\sum\limits_{N=1}^{\infty}\frac{1}{N^{2\delta-1}}\leq \int^{\infty}_{1} x^{1-2\delta} \mathrm{d}x
=\frac{1}{2\delta-2}<\infty,$$
from which we derive
$$
\sum\limits_{N=1}^{\infty}\mu^{\infty}(\Theta_N)\ll \sum\limits_{N=1}^{\infty}\frac{1}{N^{2\delta-1}}<\infty.
$$
Then Lemma \ref{BCLemma}  indicates that the set of all $\mathcal{S}$ for which $\mathcal{S}\in \Theta_N$ for infinitely many~$N$ has $\mu^{\infty}$-measure $0$. That is, with probability $1$ we have $\mathcal{S}\in \Theta_N$ for at most
finitely many $N$. So with probability $1$, we have
$$
\left|\lambda_{\mathcal{S}}(N)-\frac{N}{2}\right|\leq \delta\log_2(N)$$
for sufficiently large $N$ and the results follow. \qed

\section{Conclusions and final remarks}

In this work, we studied
 the average behavior for fixed $N$ and the asymptotic behavior for varying $N$ of the $N$th $2$-adic and rational complexities of binary sequences.

 For the average behavior,
we proved the lower bound $\frac{N}{2}-1$, which is close to the conjectured value supported by experimental results for $N\le 20$ of  \cite{TQ2010}, as well as the upper bound $\frac{N}{2}+O(\log(N))$.
We also provided analog bounds for the expected value of the rational complexity.  These results essentially answer the questions stated in \cite{GK2012,K2014,TQ2010}.

 For the asymptotic behavior,
we proved that the $2$-adic complexity profile of a random binary sequence grows approximately as $\frac{N}{2}$ with deviation of order of magnitude at most $\log(N)$ with probability $1$. This answers a question stated in \cite[p.\ 329]{GK2012}.
It also  characterizes the set of accumulation points of normalized $2$-adic complexity  $\frac{\lambda_{\mathcal{S}}(N)}{N}$ studied in \cite[Sect.\ 18.5.2, p.\ 397]{GK2012}, which
is $\{\frac{1}{2}\}$.

For periodic sequences of period $T$ the expected value of the $2$-adic complexity is
$T+O(\log(T))$, see \cite[Corollary~18.2.2]{GK2012} and \cite{HF2008}.

\section*{Acknowledgments}

The work was written during a pleasant visit of Z. Chen to Linz. He wishes to thank the Chinese Scholarship Council for financial support and  RICAM of the Austrian Academy of Sciences for hospitality.

The work was supported in part by the National Science Foundation of China under Grant No.62372256,  by the Fujian Provincial Natural Science Foundation of China under Grant No.\ 2023J01996, and by the Key Technological Innovation and Industrialization Projects in Fujian Province under Grant No.\ 2024XQ020.

\end{document}